\documentclass{article}
\usepackage{graphicx} 

\usepackage{amsmath}
\usepackage{amsfonts}

\usepackage{appendix}

\title{Expected Natural Density of Countable Sets after Infinitely Iterated de Finetti Lotteries, Computed via Matrix Decomposition}
\author{Julio Cesar Enciso-Alva}
\date{September 2024}

\newtheorem{definition}{Definition}

\newcommand{\N}{\mathbb{N}}
\newcommand{\R}{\mathbb{R}}
\newcommand{\aabs}[1]{\left| #1 \right|}
\newcommand{\card}[1]{\left| #1 \right|}
\newcommand{\sset}[1]{\left\{ #1 \right\}}
\newcommand{\ppar}[1]{\left( #1 \right)}
\newcommand{\spar}[1]{\left[ #1 \right]}

\newcommand{\floor}[1]{\left\lfloor #1 \right\rfloor}

\begin{document}
\numberwithin{equation}{section}

\maketitle

\begin{abstract}
    Consider a fair lottery over the natural numbers in which the selected number is removed. 
    This lottery is iterated countably infinite times, with a known ratio of iterations to natural numbers. Removed numbers are not replaced.
    The natural numbers are partitioned into two sets with a given ratio of elements, which is tracked along each iteration of the lottery.
    
    Hess and Polisetty considered and investigated such a process and reported the expected values of the densities for some particular cases.
    
    In this work, we provide a novel framework for computing these expected densities using infinite matrices. The results presented in this work generalize previous results. 
\end{abstract}


\section{Introduction}

A fair lottery can be defined informally as randomly selecting a ticket from a set of tickets so that all tickets are equally likely to be selected. 
Real-world lotteries work with finite sets of tickets and are typically associated with randomized decisions --including monetary prizes. 
These finite lotteries are often used to illustrate probability concepts such as equiprobability or sampling with and without substitution.

Wenmackers and Horsten \cite{fair_infinite_lottery} describe some generalizations of fair lotteries with infinite tickets in great detail.
They propose that a fair lottery should keep the following properties:
\begin{enumerate}
    \item Each single ticket shouldn't have a higher probability of being selected than any other ticket.
    \item Any individual ticket should be able to be selected.
    \item The probability of selecting a group of tickets should be equal to the sum of the probability of choosing each individual ticket.
    \item The labeling of the tickets is independent of the outcome.
\end{enumerate}

The case of a fair lottery with one ticket per natural number, $\N$, is known as the de Finetti Lottery in honor of Bruno de Finetti.
An analysis of the possibilities of modeling the de Finetti Lottery within the framework of standard Probability Theory is beyond the scope of this work, and the interested reader should refer to the paper of Wenmackers and Horsten \cite{fair_infinite_lottery}.

Hess and Polisety \cite{hess2023} proposed a variation of the de Finetti Lottery in which the tickets corresponding to odd numbers are selected with replacement, while the tickets from even numbers are selected without replacement;
this process is then iterated once for each natural number, establishing a fair lottery with the tickets remaining from the previous iteration.
This process will be referred to as the Infinitely-Iterated de Finetti Lottery (IIFL) in this text for ease of notation.

Hess and Polisety \cite{hess2023} investigated the ratio of tickets from odd numbers with respect to the remaining tickets after the iterations were performed.

In this work, I propose the concept of Finitely-Iterated de Finetti Lottery (FIFL) as a framework to study quantities related to IIFL.
By using FIFL, I was able to replicate and generalize the findings of Hess and Polisety \cite{hess2023}.


\section{Formal definitions}

This text aims to have an intuitive notion of `fair lotteries' consistent with that proposed by Wenmackers and Horsten \cite{fair_infinite_lottery}, which was described briefly in the previous section.
The discussion is limited to subsets of $\N$ and is fully contained within the standard analysis framework.


In the paper by Hess and Polisety \cite{hess2023}, a type of iterative lottery is defined, referred to in this text as an Infinitely-Iterated de Finetti Lottery (IIFL).
The ITFL starts with one ticket per each natural number, distinguishing if the number is even or odd.
On each iteration, a fair lottery is constructed with the available tickets; tickets corresponding to even numbers are selected without replacement, and tickets from odd numbers are chosen with replacement.
This process is iterated one time per each natural number.
After all the iterations were performed, Hess and Polisety investigated the expected value for the ratio (even tickets)/(natural numbers).

With the aim of self-containment, the intuitive concept of a ratio between sets is replaced with natural density; this definition is equivalent to
that used by Hess and Polisety \cite{hess2023} for
subsets of $\N$ and is compatible with the concept IIFL.

\begin{definition}[Natural density]
Let $A \subseteq \N$ be a countable set.
If the following limit is well-defined, we call it the \textbf{natural density} of $A$,
\begin{equation}
    m\ppar{A} = 
    \lim_{N\rightarrow \infty} \frac{1}{N} \card{ \sset{a\in A; a\leq N} }.
\end{equation}
\end{definition}

In the context of lotteries, this definition fails property (4)), proposed by Wenmackers and Horsten, since it doesn't guarantee that the density of any set will be kept after a permutation of labels.
This definition is used to make it compatible with the work of Hess and Polisety.

The proposed definition for the finite version of an IIFL, the Finitely Iterated de Finetti Lottery, is conceptually compatible and shares the limitation with respect to property (4).

\begin{definition}[Finitely-Iterated Lottery]
Let $A \subseteq \N$ be a countable set with $m\ppar{A} = p \in [0,1]$, and let $\alpha\geq 0$ and $N\in \N$ be parameters. 

Construct the sets $A_0 = \sset{a\in A; a\leq N}$, $B = \sset{1, 2, \dots, N} - A_0$.
The iterative lottery consists of the following steps, iterated over $n$ with $1\leq n\leq \alpha N$:
\begin{enumerate}
    \item Select randomly $x \in A_n \cup B$, with all elements being equiprobable.
    \item Construct $A_{n+1}$ as follows
    \begin{equation}
        A_{n+1} = \begin{cases}
            A_n-\sset{x}, &\text{if } x\in A_n, \\
            A_n, &\text{otherwise}.
        \end{cases}
    \end{equation}
    \item Repeat until $n< \alpha N$.
\end{enumerate}
\end{definition}

With this definition at hand, the idea of iterating the lottery `as many times as natural numbers' is formally equivalent to observing the behavior of $A_{N}$ as $N \rightarrow \infty$ with $p=\frac{1}{2}$ and $\alpha=1$.
The parameter $\alpha$ is tied to a follow-up question about performing either more or fewer iterations.

Since our interest in the lotteries is on the natural density of the outcomes instead of the actual members of the set, it is then convenient to redefine the finite de Finetti lottery regarding cardinalities.

\begin{definition}[Finitely-Iterated de Finetti Lottery (FIFL)]
Let $N\in \N$, $\pi \in [0,1]$, and $\beta \in [0, \infty)$ parameters.

The iterative lottery's result is defined by the following sequence
\begin{align}
    S\ppar{0; N, \pi} &= \floor{\pi N}, \\
    S\ppar{n+1; N, \pi} &= S\ppar{n; N, \pi} - \text{Bernoulli}\ppar{\frac{N}{N+S\ppar{n; N, \pi}}},
\end{align}
with $1\leq n\leq \beta N$.
\end{definition}

For the definition of a FIFL, $\floor{\bullet}$ is the floor function, defined as
\begin{equation}
    \floor{x} = \max\sset{k\in \N; k \leq x}
\end{equation}
and $X\sim $Bernoulli$\ppar{p}$ denotes a random variable whose probability density function is derived from the following relation
\begin{equation}
    \text{Pr}\ppar{X=x} = \begin{cases}
        p, &\text{for } x=1,\\
        1-p, &\text{for } x=0, \\
        0, &\text{otherwise}
    \end{cases}
\end{equation}

The parameters $\pi = \frac{p}{1-p}$ and $\beta = \ppar{1+\pi} \alpha$ are introduced for ease of numerical implementation. 
Their definition follows from keeping the following relations
\begin{align}
    p &=
    \frac{\card{A_0}}{\card{A_0 \cup B}} = \frac{S\ppar{0; N, \pi}}{S\ppar{0; N, \pi}+N}
    = \frac{\pi N}{\pi N + N},
    \\
    \alpha &=
    \frac{\text{iterations}}{\card{A_0 \cup B}} =
    \frac{\text{iterations}}{S\ppar{0; N, \pi}+N}
    = \frac{\beta N}{\pi N + N}.
\end{align}

Notice that the number of odd tickets, $S$, is a Markov Process since the result of the current iteration fully determines the next iteration.

Recall that we are interested in the density of odd tickets, $\bar{S}$, which can be easily computed given $S$ as
\begin{equation}
    \bar{S}\ppar{n; N, \pi} =
    \frac{S\ppar{n; N, \pi}}{S\ppar{n; N, \pi} + N}
\end{equation}

Furthermore, we are interested in the expected value of $\bar{S}$, which is denoted as follows
\begin{equation}
    m_N\ppar{\pi, \beta} = E\ppar{\bar{S}\ppar{\beta N; N, \pi }}
\end{equation}

This paper does not explore $S$ or its convergence as $N \rightarrow \infty$, but only the behavior of $m_N$. 
In particular, we consider the following quantity whenever it is well-defined
\begin{equation}
    \mu\ppar{\pi, \beta} = \lim_{N\rightarrow \infty} m_N\ppar{\pi, \beta}
\end{equation}

For readability, we define $\mu^*$ and $m_N^*$ as versions of $\mu$ and $m_N$ using the original variables $p$ and $\alpha$. 
\begin{align}
    \mu^*\ppar{p, \alpha} &= \mu\ppar{\frac{p}{1-p}, \frac{1}{1-p} \alpha}
    \\
    m_N^*\ppar{p, \alpha} &= m_N\ppar{\frac{p}{1-p}, \frac{1}{1-p} \alpha}
\end{align}
Recall that $p$ is the natural density of the set to be use on the lottery, while $\alpha$ is the ratio of iterations/density of $\N$.

Within this framework,
the findings reported by Hess and Polisety \cite{hess2023} can be described as the following
\begin{align}
\mu^*\ppar{\frac{1}{2}, 1}
    &=
    \frac{W\ppar{e^{-1}}}{1+W\ppar{e^{-1}}} \approx 0.2178
    \\
\mu\ppar{\frac{1}{2}, \frac{1}{2}}
    &=
    \frac{W\ppar{1}}{1+W\ppar{1}} \approx 0.3619
\end{align}


The proposed result in this paper is that
\begin{align}
    {\mu^*}\ppar{p, \alpha} &= \frac{W\ppar{h\ppar{p, \alpha}}}{1+W\ppar{h\ppar{p, \alpha}}}
    \\
    h\ppar{p, \alpha} &=
    \frac{p}{1-p} \exp{\ppar{\frac{p-\alpha}{1-p}}}
\end{align}
subject to the following sufficient condition
\begin{align}
    \alpha > 1 + \ppar{1-p} \ln{\ppar{\frac{p}{1-p}}} 
\end{align}

\section{Model as a Markov Process}

As described in the previous section, $S(\bullet; N, \pi)$ is dclearly a Markov Process.
The states of this process can be mapped to the possible values of $m_N$, the quantity of interest, as
\begin{equation}
    \bar{S}(t; N, \pi) \in \sset{\frac{k}{N+k}; k = 0, 1, 2, \dots}
\end{equation}
for $t\in \N$. 
Since there are countably infinite states for this process, we can define a state vector $Z_N \in \R^{\N\times 1}$ as
\begin{equation}
    \spar{Z_N(t)}(k) = Pr\ppar{\bar{S}(t; N, \pi) = \frac{k}{N+k}}, \text{ for } k=0,1, 2, \dots
\end{equation}
with the initial condition
\begin{equation}
    \spar{Z_N(0)}(k) = \spar{e_{\pi N}}(k) = \begin{cases}
        1, &\text{if } k = \max\sset{k\in \N; k \leq \pi N} \\
        0, &\text{otherwise.}
    \end{cases}
\end{equation}
with $e_\tau$ the $\tau$-th canonical vector.

The following condition gives the evolution of this system
\begin{equation}
    Pr\ppar{\bar{S}(t+1; N, \pi) = x \;\middle|\; \bar{S}(t; N, \pi) = \frac{k}{N+k}} = 
    \begin{cases}
        \frac{N}{N+k}, &\text{for } x= \frac{k}{N+k} \\
        \frac{k}{N+k}, &\text{for } x= \frac{k-1}{N+k-1} \\
        0, &\text{otherwise.}
    \end{cases}
    \label{eq:evolution}
\end{equation}
which can be encoded using a multiplication of the state vector, $Z_N$, with a transition matrix, $M_N\in \R^{\N\times \N}$, defined as
\begin{equation}
    \spar{M_N}(j,k) = 
    \begin{cases}
        \frac{N}{N+k}, &\text{if } j= k, k>0 \\
        \frac{k}{N+k}, &\text{if } j= k-1, k>0 \\
        1, &\text{if } j= k=0 \\
        0, &\text{otherwise.}
    \end{cases}
\end{equation}
Thus, the evolution equation \eqref{eq:evolution} can be rewritten as
\begin{equation}
    Z_N(t) = M_N Z_N(t-1) = \spar{M_N}^t Z_N(0).
\end{equation}

The original goal, studying the behavior of $m_N$, is connected by writing
\begin{align}
    m_N(\pi, \beta) 
    &= 
    E\ppar{\bar{S}\ppar{\beta N; N, \pi }}
    =
    \sum_{k=0}^{\infty} \frac{k}{N+k} \spar{Z_N(\beta N)}(k)
    =
    \nu^T Z_N(\beta N)
    \label{eq:weighted}
\end{align}
where $\nu_N \in \R^{\N \times 1}$ is a vector given by
\begin{equation}
    \spar{\nu_N}(k) = \frac{k}{N+k}, \text{ for } k=0,1, 2, \dots
\end{equation}

Notice that, in every finite case, the infinite sum in \eqref{eq:weighted} can be truncated at $k\leq \beta N$ since $Z_N$ is zero for those values.
However, it is convenient to keep it like that since it allows us to consider all possible values of $\pi$ by changing only the initial conditions.
This becomes clear by rewriting $m_N$ as
\begin{equation}
    m_N(\pi, \beta) = \nu_N^T \spar{M_N}^{\beta N}\, e_{\pi N}.
    \label{eq:mat_mult_base}
\end{equation}


\section{Results}
\label{sec:result}

Our strategy, whose details are detailed in the appendix, is to compute a sort of eigendecomposition of $M_N$ and use it to compute $m_N$ efficiently.
The existence of such objects is guaranteed since $M_N$ is lower diagonal, and thus its eigenvalues are trivially found as
\begin{equation}
    \lambda^{(N)}_k = \frac{N}{N+k}, \text{ for } k=0,1, 2, \dots
\end{equation}

The $k$-th eigenvalue of $M_N$, which we now refer as $V_k^{(k)}$, is defined by the following property
\begin{equation}
    M_N V_k^{(N)} = \lambda^{(N)}_k V_k^{(N)}
\end{equation}

As it is shown in the appendix \ref{ap:eigenvalues}, those eigenvalues are given by
\begin{equation}
    \spar{V_k^{(N)}}(n) =
    (-1)^{k-n} \ppar{\frac{N+k}{N}}^{k-n} \cdot \frac{N+n}{N+k} \binom{k}{n}
\end{equation}

With this result at hand, the eigendecomposition of $M_N$ is given by
\begin{align}
    M_N &= \spar{V_N}^{-1} \Lambda_N V_N 
    \label{eq:eigendecomp}
    \\
    V_N &= \begin{bmatrix}
        V_0^{(N)}, & V_1^{(N)}, & \dots
    \end{bmatrix} \\
    \Lambda_N &= \text{diag}\ppar{\lambda^{(N)}_0, \lambda^{(N)}_1, \dots}
\end{align}

As it is proven in the appendix \ref{ap:inverse}, the inverse of $V_N$ exists and is given by
\begin{align}
    \spar{V^{-1}}(n,k) &=
    \ppar{\frac{N+n}{N}}^{k-n} \binom{k}{n}
\end{align}

The objective of defining and computing the eigendecomposition of $M_N$ is to compute efficiently $m_N$.
In particular, from equations \eqref{eq:mat_mult_base} and \eqref{eq:eigendecomp} we have
\begin{equation}
    m_N(\pi, \beta) = \nu_N^T 
    \spar{V_N}^{-1} \spar{\Lambda_N}^{\beta N} V_N\,
    e_{\pi N}
\end{equation}

As it is proven in the appendix \ref{ap:product}, the result of this computation is given by
\begin{equation}
m_N(\pi, \beta) =
    \sum_{k=0}^{\pi N}
    (-1)^{k+1} 
    \binom{\pi N}{k}
    \ppar{\frac{k}{N+k}}^k
    \ppar{\frac{N+k}{N}}^{\spar{\pi-\beta}k N -1}
    \label{eq:eq2}
\end{equation}

As it is proven in the appendix \ref{ap:coeff}, the coefficients on equation \eqref{eq:eq2} converge pointwise with $N\rightarrow \infty$ to the following expression
\begin{equation}
    \mu(\alpha, \beta) =
    \sum_{k=0}^{\infty}
    (-1)^{k+1} \frac{k^k}{k!} 
    \ppar{\pi
    e^{\spar{\pi-\beta}} }^k
    \label{eq:mu_sum}
\end{equation}

As it is proven in the appendix \ref{ap:convergence},
the series on \eqref{eq:mu_sum} converges under the following condition
\begin{equation}
    \beta > \pi + 1 + \ln{\pi}
\end{equation}

A more compact expression for \eqref{eq:mu_sum} can be obtained using the following identity, derived in appendix \ref{ap:identities},
\begin{equation}
    \frac{W\ppar{z}}{1+W\ppar{z}}
    =
    \sum_{k=0}^{\infty}
    (-1)^{k+1} \frac{k^k}{k!} 
    z^k
\end{equation}
resulting in the following 
\begin{align}
    \mu\ppar{\pi, \beta} &= \frac{W\ppar{g\ppar{\pi, \beta}}}{1+W\ppar{g\ppar{\pi, \beta}}}
    \\
    g\ppar{\pi, \beta} &=
    \pi \exp{\ppar{\pi - \beta}}
\end{align}

It is worth to rewrite this result in terms of $p$ and $\alpha$, the original variables.
\begin{align}
    {\mu^*}\ppar{p, \alpha} &= \frac{W\ppar{h\ppar{p, \alpha}}}{1+W\ppar{h\ppar{p, \alpha}}}
    \label{eq:final1}
    \\
    h\ppar{p, \alpha} &=
    \frac{p}{1-p} \exp{\ppar{\frac{p-\alpha}{1-p}}}
    \label{eq:final2}
\end{align}
subject to the following
\begin{align}
    \alpha > 1 + \ppar{1-p} \ln{\ppar{\frac{p}{1-p}}} 
    \label{eq:conv_og}
\end{align}

\section{Discussion}

\begin{figure}
    \centering
    \includegraphics[width=\linewidth]{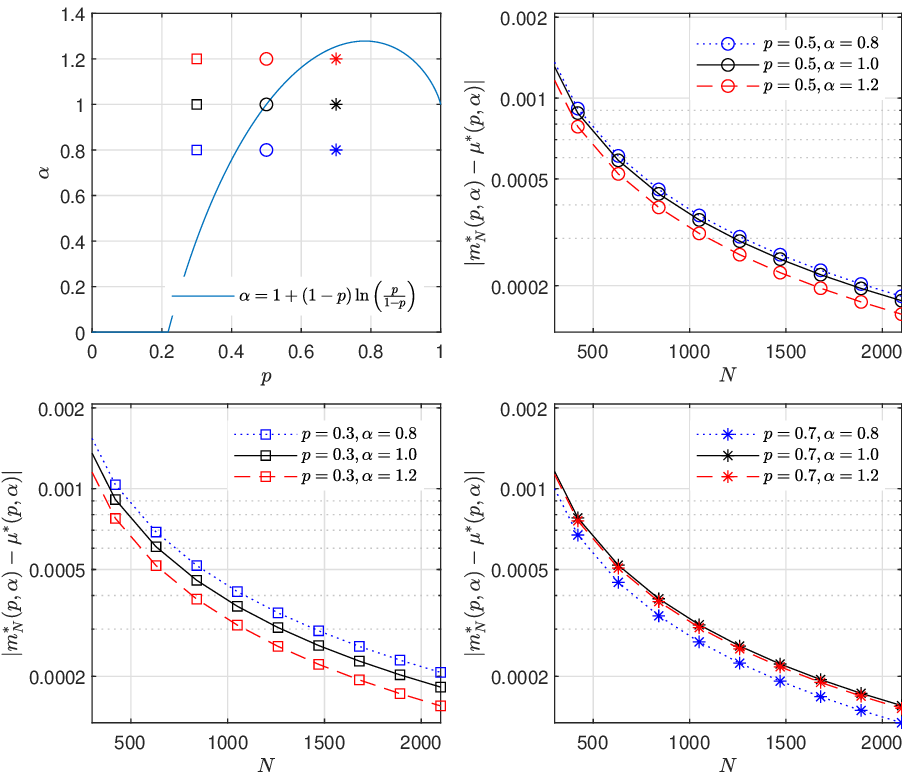}
    \caption{Numerical evidence of the convergence of $m^*_N$ as $N\rightarrow \infty$ for different values of $p$ and $\alpha$. Notice that this apparently contradicts the condition \eqref{eq:conv_og}; more details are provided in the text.}
    \label{fig:onlyfig}
\end{figure}

The formula on equations \eqref{eq:final1} and \eqref{eq:final2} is consistent with the findings by Hess and Polisetty, and it constitutes a slight generalization.

%
The processes described in section \ref{sec:result} were implemented in Matlab in order to perform the computations for large values of $N$ (more than 2,000).
Although the condition on equation \eqref{eq:conv_og} was found sufficient to guarantee the convergence of the computed quantities, the numerical computations show that the condition may not be necessary.
For instance, in Figure \ref{fig:onlyfig} are displayed some cases that should not converge but do, in fact, converge.

One possible explanation for this apparent contradiction is that the convergence condition on \eqref{eq:conv_og} was determined after taking the convergence of the coefficient as $N\rightarrow \infty$.
%
%
It can be easily verified
that the series that defines $\mu^*$ is heavy-tailed, whereas each one of $m^*_N$ is not.

Thus, we formally state that the established condition is sufficient but also provide informal evidence that it is not necessary.
Developing a sufficient and necessary condition was deemed out of the scope of this work.


Under more careful handling, this technique can be used for several variations of the same problem: duplicating the selected elements instead of removing them, removing multiples of the selected elements, multi-staged removal after repeated selection, etc.

On the other hand, this technique lacks the granularity to manage individual elements. 
%
%
%
This is consistent with violating the property (4) of fairness: independence of label permutations.



\bibliographystyle{plain} 
\bibliography{finetti} 


\newpage

\appendix 

\section{Eigenvalues of Transition Matrix}
\label{ap:eigenvalues}

Recall that the transition matrix, $M_N\in \R^{\N\times \N}$, is defined as
\begin{equation}
    \spar{M_N}(j,k) = 
    \begin{cases}
        \frac{N}{N+k}, &\text{if } j= k, k>0 \\
        \frac{k}{N+k}, &\text{if } j= k-1, k>0 \\
        1, &\text{if } j= k=0 \\
        0, &\text{otherwise.}
    \end{cases}
\end{equation}

The eigenvalues of $M_N$ are given by
\begin{equation}
    \lambda^{(N)}_k = \frac{N}{N+k}, \text{ for } k=0,1, 2, \dots
\end{equation}

The $k$-th eigenvalue of $M_N$, which we now refer as $V_k^{(k)}$, is defined by the following property
\begin{equation}
    M_N V_k^{(N)} - \lambda^{(N)}_k V_k^{(N)} = 0
    \label{a:01}
\end{equation}

The $j$-th column of \eqref{a:01} is given by
\begin{equation}
    \ppar{\frac{N}{N+j}-\frac{N}{N+k} } V_k^{(N)}(j) + \ppar{\frac{j+1}{N+j+1}}V_k^{(N)}(j+1) = 0
    \label{a:induction}
\end{equation}

One particular case is the $k$-th column,
\begin{equation}
    \ppar{\frac{k+1}{N+k+1}}V_k^{(N)}(k+1) = 0
\end{equation}
from which it follows that $V_k^{(N)}(k+1) = 0$. Furthermore, by induction over \eqref{a:induction} this implies that $V_k^{(N)}(j) = 0$ for $j>k$.

For $j\leq k$, equation \eqref{a:induction} can be rewritten as
\begin{align}
    V_k^{(N)}(j) 
    &= -
    \frac{\ppar{\frac{j+1}{N+j+1}}}{\ppar{\frac{N}{N+j}-\frac{N}{N+k} }}
    V_k^{(N)}(j+1)
    \nonumber \\
    &=
    -
    \ppar{\frac{N+k}{N}} \ppar{\frac{N+j}{N+j-1}} \ppar{\frac{j+1}{k-j}} V_k^{(N)}(j+1)
\end{align}

Without loss of generality, we can define $\spar{V_k^{(N)}}(k) = 1$, and then the other entries are given by
\begin{equation}
    \spar{V_k^{(N)}}(j) =
    \begin{cases}
    (-1)^{k-j} \ppar{\frac{N+k}{N}}^{k-j} \cdot \frac{N+j}{N+k} \binom{k}{j},
    &\text{for } 0 \leq j\leq k\\
    0, &\text{otherwise.}
    \end{cases}
\end{equation}

Once the eigenvalues are found, it is possible to ensemble a matrix $V_N$ whose columns are the eigenvalues of $M_N$.
\begin{equation}
    V_N = \begin{bmatrix}
        V_0^{(N)}, & V_1^{(N)}, & \dots
    \end{bmatrix}
\end{equation}

\section{Inverse of Eigenvalues of Transition Matrix}
\label{ap:inverse}

The purpose of computing the eigendecomposition of $M_N$ is to use the identity
\begin{equation}
    \spar{M_N}^t = V_N \spar{\Lambda_N}^t \spar{V_N}^{-1}
\end{equation}
it is sufficient that $\spar{V_N}^{-1}$ is a left inverse of $V_N$. 

Our strategy for computing this operator is using back-substitution.
For ease of notation, consider the inverse of $V_N$ as a collection of row vectors
\begin{equation}
    \spar{V_N}^{-1} = \begin{bmatrix}
        U_0^{(N)} \\ U_1^{(N)} \\ \vdots
    \end{bmatrix}
    \label{b:00}
\end{equation} 

With this notation, to find a left inverse of $V_N$ is to solve to solve the following system
\begin{equation}
    U_k^{(N)}\, V_N = e_k
    \label{b:01}
\end{equation}
for $k = 0, 1, \dots$, where $e_k\in \R^\N$ the $k$-th canonical vector given by
\begin{equation}
    \spar{e_k}(j) = \begin{cases}
        1, &\text{if } j=k, \\
        0, &\text{otherwise}
    \end{cases}
\end{equation}

Since $V_N$ is a lower diagonal matrix, we may assume that $\spar{V_N}^{-1}$ is too, and so we define 
\begin{equation}
    {U_k^{(N)}}(j) = 0, \text{ for } j<k
\end{equation}
which trivially satisfies the $j$-th column of equation \eqref{b:01} for $j<k$.
On a similar way, the $k$-th column of equation \eqref{b:01} is satisfied if we select
\begin{equation}
    {U_k^{(N)}}(k) = 1
\end{equation}

In general, the $j$-th column of equation \eqref{b:01} can be written as
\begin{equation}
    \sum_{m=0}^\infty {U_k^{(N)}}(m)\, V(m,j)  = 
    \begin{cases}
        1, &\text{if } j=k \\
        0, &\text{otherwise}
    \end{cases}
    \label{b:02}
\end{equation}

The back-substitution step follows from equation \eqref{b:02};
we shall use the expression
\begin{equation}
    {U_k^{(N)}}(j) = 
    -\sum_{m=k}^{j-1} {U_k^{(N)}}(m)\, V(m,j)
\end{equation}
to compute the entries ${U_k^{(N)}}$ using entries previously computed, starting with ${U_k^{(N)}}(k) = 1$.
After some ad hoc induction, not shown here, we find the following expression
\begin{align}
U_k^{(N)}(j)
    &=
    \ppar{\frac{N+k}{N}}^{j-k} \binom{j}{k}
    \label{eq:b03}
\end{align}

Now we prove that $U_k^{(N)}$ from \eqref{eq:b03} satisfies the equation \eqref{b:02} for $k<j$.
%
\begin{align*}
    \sum_{m=k}^j {U_k^{(N)}}(m)\, V(m,j)
    &=
    \sum_{m=k}^j
    \spar{ \ppar{\frac{N+k}{N}}^{m-k} \binom{m}{k} }
    \\
    &\phantom{=}
    \phantom{\sum_{m=k}^j} \cdot
    \spar{ (-1)^{j-m} \ppar{\frac{N+j}{N}}^{j-m} \frac{N+m}{N+j} \binom{j}{m} }
    \\
    &=
    \frac{(N+j)^j}{(N+k)^k} \cdot \frac{1}{N^{j-k}(N+j)} \cdot \frac{j!}{k!}
    \\
    &\phantom{=}
    \phantom{\sum_{m=k}^j} \cdot
    \sum_{m=k}^j
    (-1)^{j-m} \ppar{\frac{N+k}{N+j}}^m \frac{(N+m)}{(j-m)!\, (m-k)!}
    \\
    &=
    \frac{(N+j)^{j-k}}{N^{j-k}(N+j)} \cdot \frac{j!}{k!\, (j-k)!} \cdot (-1)^{j-k}
    \\
    &\phantom{=}
    \phantom{\sum_{m=k}^j} \cdot
    \sum_{m=k}^{j}
    (-1)^{m-k} \ppar{\frac{N+k}{N+j}}^{m-k} 
    (N+m) \binom{j-k}{m-k}
    \\
    &=
    \ppar{K_{j,k} }
    \sum_{\tau=0}^{j-k}
    \binom{j-k}{\tau}
    \ppar{-\frac{N+k}{N+j}}^{\tau} 
    (N+k + \tau) 
    \\
    &=
    \ppar{K_{j,k} } 
    \sum_{\tau=0}^{j-k}
    \binom{j-k}{\tau}
    \ppar{-\frac{N+k}{N+j}}^{\tau} 
    \tau
    \\
    &\phantom{=}
    \phantom{\sum_{m=k}^j}
    +
    \ppar{K_{j,k} } \ppar{N+k}
    \sum_{\tau=0}^{j-k}
    \binom{j-k}{\tau}
    \ppar{-\frac{N+k}{N+j}}^{\tau} 
    \\
    &=
    \ppar{K_{j,k} } 
    \sum_{\tau=1}^{j-k}
    (j-k)
    \binom{j-k-1}{\tau-1}
    \ppar{-\frac{N+k}{N+j}}^{\tau} 
    \\
    &\phantom{=}
    \phantom{\sum_{m=k}^j}
    +
    \ppar{K_{j,k} } \ppar{N+k}
    \ppar{1-\frac{N+k}{N+j}}^{j-k}
    \\
    &=
    \ppar{K_{j,k} } 
    (j-k)
    \ppar{-\frac{N+k}{N+j}}
    \ppar{1-\frac{N+k}{N+j}}^{j-k-1}
    \\
    &\phantom{=}
    \phantom{\sum_{m=k}^j}
    +
    \ppar{K_{j,k} } \ppar{N+k}
    \ppar{1-\frac{N+k}{N+j}}^{j-k}
    \\
    &=
    \ppar{K_{j,k} } \ppar{1-\frac{N+k}{N+j}}^{j-k-1}
    \\
    &\phantom{=}
    \phantom{\sum_{m=k}^j}
    \cdot
    \spar{(j-k)\ppar{-\frac{N+k}{N+j}}+(N+k)\ppar{1-\frac{N+k}{N+j}}}
    \\
    &=
    \ppar{K_{j,k} } \ppar{\frac{j-k}{N+j}}^{j-k-1} 
    \frac{N+k}{N+j}
    \spar{-(j-k) + (j-k)}
    \\ &= 0
\end{align*}

The variable ${\ppar{K_{j,k} }}$ is used to shorten terms that depend only on $j$ and $k$.
After this computation, we conclude that, for $k<j$,
\begin{equation}
    \sum_{m=k}^j {U_k^{(N)}}(m)\, V(m,j) = 0
\end{equation}
and so ${U_k^{(N)}}$ satisfy equation \eqref{b:01}.
Then $\spar{V_N}^{-1}$, which is constructed as in equation \eqref{b:00}, is a left inverse of $V_N$.

\section{Expected Value of Matrix Product}
\label{ap:product}

Recall that the goal of using the eigendecomposition of the transition matrix is to compute the expected density efficiently.
In other words, the goal is to perform the following computation
\begin{equation}
    m_N(\pi, \beta) = 
    {\nu_N^T 
    V_N}\, \spar{\Lambda_N}^{\beta N} \spar{V_N}^{-1}\,e_{\pi N}
\end{equation}

%
To ease the computation further, it is carried by parts, starting with the product $\nu_N^T V_T$,
\begin{align*}
    \spar{\nu_N^T V_N}(m) &=
    \sum_{k=0}^\infty \spar{\nu_N}(k) \cdot V_N(k,m) \\
    &=
    \sum_{k=0}^m
    \ppar{\frac{k}{N+k}}
    \ppar{(-1)^{m-k} \ppar{\frac{N+m}{N}}^{m-k} \cdot \frac{N+k}{N+m} \binom{m}{k}}
    \\
    &=
    \frac{1}{N+m}
    \sum_{k=0}^m
    \binom{m}{k}
    \ppar{- \frac{N+m}{N}}^{m-k} k
    \\
    &=
    \frac{m}{N+m}
    \sum_{k=1}^m
    \binom{m-1}{k-1}
    \ppar{- \frac{N+m}{N}}^{m-k}
    \\
    &=
    \frac{m}{N+m}
    \sum_{k'=0}^{m-1}
    \binom{m-1}{k'}
    \ppar{- \frac{N+m}{N}}^{m-1-k'}
    \\
    &=
    \frac{m}{N+m}
    \ppar{1-\frac{N+m}{N}}^{m-1}
    \\
    &=
    \frac{m}{N+m}
    \ppar{-\frac{m}{N}}^m \ppar{-\frac{N}{m}}
    \\
    &=
    (-1)^{m-1} \ppar{\frac{m}{N}}^m \frac{N}{N+m}
\end{align*}

The matrix $\spar{\Lambda_N}^{\beta N}$ is trivially a diagonal matrix whose entries are given by
\begin{equation}
    \spar{\ppar{\Lambda_N}^{\beta N}}(m,m)
    =
    \ppar{\frac{N}{N+m}}^{\beta N}
\end{equation}

Thus the product $\nu_N^T V_T \spar{\Lambda_N}^{\beta N}$ is almost trivial
\begin{equation}
    \spar{\nu_N^T V_T \spar{\Lambda_N}^{\beta N}}(m)
    =
    (-1)^{m-1} \ppar{\frac{m}{N}}^m
    \ppar{\frac{N}{N+m}}^{\beta N+1}
\end{equation}

The product $\spar{V_N}^{-1}\,e_{\pi N}$ is equivalent to selecting the $\pi N$-th column of $\spar{V_N}^{-1}$, which is given by
\begin{align}
\spar{V_N}^{-1}(m, \pi N)
    &=
    \begin{cases}
        \ppar{\frac{N+m}{N}}^{\pi N-m} \binom{\pi N}{m},
        &\text{if } m<\pi N \\
        0, &\text{otherwise}
    \end{cases}
\end{align}


Finally, the quantity of interest is computed
\begin{align*}
    m_N(\pi, \beta) 
    &=
    \ppar{\nu_N^T V_N \spar{\Lambda_N}^{\beta N}} \ppar{\spar{V_N}^{-1}\,e_{\pi N}}
    \\
    &=
    \sum_{m=0}^\infty 
    \spar{\nu_N^T V_T \spar{\Lambda_N}^{\beta N}}(m) \cdot
    \spar{V_N}^{-1}(m, \pi N)
    \\
    &=
    \sum_{m=0}^{\pi N}
    \ppar{(-1)^{m-1} \ppar{\frac{m}{N}}^m
    \ppar{\frac{N}{N+m}}^{\beta N+1}}
    \\
    &\phantom{=}
    \phantom{\sum_{m=0}^{\pi N}} \cdot
    \ppar{\ppar{\frac{N+m}{N}}^{\pi N-m} \binom{\pi N}{m}}
    \\
    &=
    \sum_{m=0}^{\pi N}
    (-1)^{m-1} \binom{\pi N}{m}
    \ppar{\frac{m}{N}}^m \ppar{\frac{N+m}{N}}^{-m}
    \\
    &\phantom{=}
    \phantom{\sum_{m=0}^{\pi N}} \cdot
    \ppar{\frac{N}{N+m}}^{\beta N+1}
    \ppar{\frac{N+m}{N}}^{\pi N}
    \\
    &=
    \sum_{k=0}^{\pi N}
    (-1)^{k+1} 
    \binom{\pi N}{k}
    \ppar{\frac{k}{N+k}}^k
    \ppar{\frac{N+k}{N}}^{\spar{\pi-\beta} N -1}
\end{align*}

\section{Pointwise Convergence of Coefficients}
\label{ap:coeff}

Now that we have the following expression for the quantity of interest,
\begin{equation}
    m_N(\pi, \beta) =
    \sum_{k=0}^{\pi N}
    (-1)^{k+1} 
    \binom{\pi N}{k}
    \ppar{\frac{k}{N+k}}^k
    \ppar{\frac{N+k}{N}}^{\spar{\pi-\beta} N -1}
\end{equation}
it is relevant to study its behavior as $N\rightarrow \infty$. 
To be precise, we are interested in the following limit
\begin{multline}
    \lim_{N \rightarrow \infty}
    (-1)^{k+1} 
    \binom{\pi N}{k}
    \ppar{\frac{k}{N+k}}^k
    \ppar{\frac{N+k}{N}}^{\spar{\pi-\beta} N -1}
    \\
    =
    (-1)^{k+1}
    \ppar{\lim_{N \rightarrow \infty}
    \binom{\pi N}{k}
    \ppar{\frac{k}{N+k}}^k}
    \ppar{\lim_{N \rightarrow \infty}
    \ppar{\frac{N+k}{N}}^{\spar{\pi-\beta} N -1}}
\end{multline}

The second limit can be dealt with easily
\begin{align*}
    \lim_{N \rightarrow \infty}
    \ppar{\frac{N+k}{N}}^{\spar{\pi-\beta} N -1}
    &=
    \lim_{N \rightarrow \infty}
    \ppar{\ppar{1+\frac{k}{N}}^N}^{\spar{\pi - \beta}}
    \ppar{\frac{N+k}{N}}
    =
    e^{\spar{\pi - \beta}k}
\end{align*}

The first limit can be approximated using a technique used on the Poisson approximation of binomial random variables. 
To be precise
\begin{align*}
    \lim_{N \rightarrow \infty}
    \binom{\pi N}{k}
    \ppar{\frac{k}{N+k}}^k
    &=
    \lim_{N \rightarrow \infty}
    \frac{\ppar{\pi N}!}{k!\, \ppar{\pi N-k}!} \cdot 
    \frac{k^k}{\ppar{N+k}^k}
    \\
    &=
    \frac{k^k}{k!}
    \lim_{N \rightarrow \infty}
    \frac{\ppar{\pi N}!}{\ppar{\pi N-k}!\, \ppar{N+k}^k}
    \\
    &=
    \frac{k^k}{k!}
    \lim_{N \rightarrow \infty}
    \frac{\ppar{\pi N-k} \cdots \ppar{\pi N} }{\ppar{N+k}^k}
    \\
    &=
    \frac{k^k}{k!}
    \lim_{N \rightarrow \infty}
    \ppar{\frac{\pi N-k}{N+k}} \cdots \ppar{\frac{\pi N}{N+k}}
    \\
    &=
    \frac{k^k}{k!}
    \pi ^k
\end{align*}

Thus, the intended limit is found
\begin{multline}
    \lim_{N \rightarrow \infty}
    (-1)^{k+1} 
    \binom{\pi N}{k}
    \ppar{\frac{k}{N+k}}^k
    \ppar{\frac{N+k}{N}}^{\spar{\pi-\beta} N -1}
    =
    (-1)^{k+1}
    \frac{k^k}{k!}
    \pi ^k
    e^{\spar{\pi - \beta}k}
\end{multline}

\section{Convergence of Density Function}
\label{ap:convergence}

the proposed density function is given by
\begin{equation}
    \mu(\alpha, \beta) =
    \sum_{k=0}^{\infty}
    (-1)^{k+1} \frac{k^k}{k!} 
    \ppar{\pi
    e^{\spar{\pi-\beta}} }^k
\end{equation}

This series's convergence is investigated using the ratio test over its coefficients.
\begin{align*}
    1 &>
    \lim_{k\rightarrow \infty}
    \frac{\aabs{
    (-1)^{\ppar{k+1}+1} \frac{\ppar{k+1}^{\ppar{k+1}}}{\ppar{k+1}!} 
    \ppar{\pi
    e^{\spar{\pi-\beta}} }^{\ppar{k+1}}
    }}{\aabs{
    (-1)^{k+1} \frac{k^k}{k!} 
    \ppar{\pi
    e^{\spar{\pi-\beta}} }^k
    }}
    \\
    &=
    \lim_{k\rightarrow \infty}
    {
    \frac{\ppar{k+1}^{\ppar{k+1}}\, k! }{{k}^{k} \ppar{k+1}!}
    {\pi e^{\spar{\pi-\beta}} }
    }
    \\
    &=
    \lim_{k\rightarrow \infty}
    \ppar{\frac{k+1}{k}}^k {\pi e^{\spar{\pi-\beta}} }
    \\
    &=
    {\pi e^{\spar{\pi-\beta}+1} }
\end{align*}

Taking logarithm, the expression is simplified to
\begin{equation}
    0 > \ln{\pi} + \spar{\pi-\beta}+1
\end{equation}
which is equivalent to
\begin{equation}
    \beta > \pi + 1 + \ln{\pi}
\end{equation}

\section{Density Function and the Lambert W Function}
\label{ap:identities}

The Lambert W Function, $W: \R_+\rightarrow \R$, is defined\footnote{The domain of the Lambert W function is constrained for the scope of this work.} by the following relation
\begin{equation}
    W(z)\, e^{W(z)} = z
\end{equation}

Using implicit differentiation over this expression leads to an identity that will be used later.
\begin{align*}
    1 &=
    \frac{d}{d z} \spar{W(z) e^{W(z)}}
    \\ &=
    W'(z)\, e^{W(z)} + W(z)\, e^{W(z)}\, W'(z)
    \\ &=
    W'(z) \ppar{1+W\ppar{z}} e^{W(z)}
    \\ &=
    W'(z) \ppar{1+W\ppar{z}} \frac{z}{W(z)}
\end{align*}

After simplification, we have the following identity
\begin{equation}
    \frac{d}{dz} W(z) 
    =
    \frac{1}{z} \cdot \frac{W(z)}{1+W(z)}
\end{equation}
Although this identity seems to come out of nowhere, it is motivated by the results obtained by Hess and Poliset \cite{hess2023}.

Now, another important identity is the Taylor expansion of the Lambert W function, given by
\begin{align}
    W(z) 
    &=
    \sum_{k=1}^\infty 
    \frac{(-k)^{k-1}}{k!} z^k
\end{align}

With both identities at hand, we can establish the following computation
\begin{align*}
    \frac{W(z)}{1+W(z)}
    &=
    z \cdot \frac{d}{dz} W(z)
    \\
    &=
    z \cdot \frac{d}{dz} \ppar{\sum_{k=1}^\infty 
    \frac{(-k)^{k-1}}{k!} z^k}
    \\
    &=
    z \cdot {\sum_{k=1}^\infty 
    \ppar{-1}^{k-1}
    \frac{k^{k-1}}{\ppar{k-1}!} z^{k-1}}
    \\
    &=
    {\sum_{k=1}^\infty 
    \ppar{-1}^{k-1}
    \frac{k^{k}}{{k}!} z^{k}}
\end{align*}

The resulting identity is
\begin{equation}
    \frac{W(z)}{1+W(z)}
    =
    {\sum_{k=1}^\infty 
    \ppar{-1}^{k-1}
    \frac{k^{k}}{{k}!} z^{k}}
\end{equation}

Notice the resemblance with the density function
\begin{equation}
    \mu(\alpha, \beta) =
    \sum_{k=0}^{\infty}
    (-1)^{k+1} \frac{k^k}{k!} 
    \ppar{\pi e^{\spar{\pi-\beta}} }^k
\end{equation}
which is why we can state the following
\begin{equation}
    \mu(\alpha, \beta) =
    \frac{W\ppar{\pi e^{\spar{\pi-\beta}}}}{1+W\ppar{\pi e^{\spar{\pi-\beta}}}}
\end{equation}


\end{document}